\documentclass[12pt,a4paper,twoside]{article}
\usepackage{amsfonts}
\usepackage{amsthm,amsmath,amssymb}

\def\trait #1 #2 #3 {\vrule width #1pt height #2pt depth #3pt}
\def\fin{
    \trait .3 5 0
    \trait 5 .3 0
    \kern-5pt
    \trait 5 5 -4.7
    \trait 0.3 5 0
\medskip}

\newtheorem{teor}{Theorem}[section]
\newtheorem{defin}[teor]{Definition}
\newtheorem{lemm}[teor]{Lemma}
\newtheorem{osse}[teor]{Remark}
\newtheorem{prop}[teor]{Proposition}
\newtheorem{defi}[teor]{Definition}
\newtheorem{coro}[teor]{Corollary}
\newtheorem{prob}[teor]{Problem}
\newtheorem{hypo}[teor]{Hypothesis}
\newcommand{\bele}{\begin{lemm}\begin{sl}}
\newcommand{\enle}{\end{sl}\end{lemm}}
\newcommand{\bedef}{\begin{defi}\begin{sl}}
\newcommand{\eddef}{\end{sl}\end{defi}}
\newcommand{\bete}{\begin{teor}\begin{sl}}
\newcommand{\ente}{\end{sl}\end{teor}}
\newcommand{\beos}{\begin{osse}\begin{sl}}
\newcommand{\eddos}{\end{sl}\end{osse}}
\newcommand{\bepr}{\begin{prop}\begin{sl}}
\newcommand{\empr}{\end{sl}\end{prop}}
\newcommand{\bepro}{\begin{prob}\begin{rm}}
\newcommand{\empro}{\end{rm}\end{prob}}
\newcommand{\bede}{\begin{defin}\begin{sl}}
\newcommand{\edde}{\end{sl}\end{defin}}
\newcommand{\beco}{\begin{coro}\begin{sl}}
\newcommand{\enco}{\end{sl}\end{coro}}
\newcommand{\behy}{\begin{hypo}\begin{sl}}
\newcommand{\enhy}{\end{sl}\end{hypo}}

\newcommand{\thspace}{\hspace{3mm}}

\newcommand{\RR}{\mathbb{R}}

\newcommand{\beeq}[1]{\begin{equation}\label{#1}}
\newcommand{\eddeq}{\end{equation}}
\newcommand{\beeqa}[1]{\begin{eqnarray}\label{#1}}
\newcommand{\eddeqa}{\end{eqnarray}}
\newcommand{\beal}[1]{\begin{align}\label{#1}}
\newcommand{\eddal}{\end{align}}
\newcommand{\bespl}[1]{\begin{split}\label{#1}}
\newcommand{\edspl}{\end{split}}
\newcommand{\bega}[1]{\begin{gather}\label{#1}}
\newcommand{\edga}{\end{gather}}
\newcommand{\beeqax}{\begin{eqnarray*}}
\newcommand{\eddeqax}{\end{eqnarray*}}
\def\qed{\ifmmode   \else \leavevmode\unskip\penalty9999 \hbox{}\nobreak\hfill
  \fi
  \quad\hbox{\hskip.5em\vrule width.4em height.6em depth.05em\hskip.1em}}
\def\endproofsym{\qed}

\def\endnobox{\def\endproofsym{}\end{proof}\def\endproofsym{\qed}}

\newcommand{\no}{\nonumber}
\newcommand{\beeqao}{\begin{eqnarray}\no}
\newcommand{\bealo}{\begin{align}\no}
\newcommand{\besplo}{\begin{split}\no}
\newcommand{\begao}{\begin{gather}\no}
\def\trait #1 #2 #3 {\vrule width #1pt height #2pt depth #3pt}
\def\fin{\hfill
    \trait .3 5 0
    \trait 5 .3 0
    \kern-5pt
    \trait 5 5 -4.7
    \trait 0.3 5 0
\medskip}

\newcommand{\duav}[1]{\langle{#1}\rangle}

\newcommand{\dt}{\partial_t}

\newcommand{\dn}{\partial_{\bf n}}

\newcommand{\itt}{\int_0^t}

\newcommand{\io}{\int_\Omega}

\newcommand{\e}{\varepsilon}
\newcommand{\ee}{^{\varepsilon}}

\newcommand{\mezzo}{{\frac{1}{2}}}

\newcommand{\fhi}{\varphi}

\newcommand{\lsc}{lower semicontinuous}

\DeclareMathOperator{\dive}{div}

\newcommand{\HUH}{H^1(0,T;H)}

\newcommand{\HUVp}{H^1(0,T;V')}
\newcommand{\CZH}{C^0([0,T];H)}
\newcommand{\CZV}{C^0([0,T];V)}

\newcommand{\LDH}{L^2(0,T;H)}
\newcommand{\LDV}{L^2(0,T;V)}
\newcommand{\LDVp}{L^2(0,T;V')}
\newcommand{\LIH}{L^\infty(0,T;H)}
\newcommand{\LIV}{L^\infty(0,T;V)}

\newcommand{\LDW}{L^2(0,T;W)}

\newcommand{\CZtH}{C^0([0,t];H)}

\newcommand{\LDtH}{L^2(0,t;H)}
\newcommand{\LDtV}{L^2(0,t;V)}
\newcommand{\LDtVp}{L^2(0,t;V')}
\newcommand{\LItH}{L^\infty(0,t;H)}

\newcommand{\betapb}{\bar{\beta_1}}
\newcommand{\betasb}{\bar{\beta_2}}
\newcommand{\xib}{\boldsymbol{\xi}}

\newcommand{\betab}{{\boldsymbol{\beta}}}

\newcommand{\teta}{\vartheta}

\let\sz s
\def\ss{^{\sigma}}
\def\s{{\sigma}}
\def\ee{^{\e}}

\def\fhi{\varphi}

\def\fine{\hfill\kern4pt \vrule height4pt depth0pt width4pt }

\numberwithin{equation}{section}
\numberwithin{equation}{section}
\begin{document}
%\date{}
%%%%%%%%%%%%%%%%%%%%%%%%%%%%%%%%%%%%%%%%%%%%%%%%%%%%%%%%%%%%%%%%%%%
\font \corr = cmitt10 scaled \magstep1
\font \correr = cmitt10 scaled \magstep1
\font \corre = cmssbx10 scaled \magstep1
\newcommand{\corrmath}[1]{\mbox{\mathversion{bold}$\mathbf{#1}$}}

\newcommand\fer[1]{\noindent {\correr \fbox{
\begin{minipage}[t]{3.8cm}
COMMENT ER :
\end{minipage}
\begin{minipage}[t]{10.5cm}
{#1}
\end{minipage}
}}}

%%%%%%%%%%%%%%%%%%%%%%%%%%%%%%%%%%%%%%%%%%%%%%%%%%%%%%%%%%%%%%%%%%%

\title{Solid liquid phase changes with different densities}

\author{\renewcommand{\thefootnote}{\!$\fnsymbol{footnote}$}
Michel Fr\'emond\\
{\sl  Dipartimento di Ingegneria Civile, Universit\`a di Roma ``Tor Vergata''}\\
 {\sl Via del Politecnico, 1, I-00133 Roma, Italy}\\
{\sl  CMLA, ENS Cachan -D\'epartement de M\'ecanique, ENSTA, Paris,}\\
{\rm E-mail:~~\tt fremond@cmla.ens-cachan.fr}\\
\and
Elisabetta Rocca\\
{\sl Dipartimento di Matematica, Universit\`a di Milano}\\
{\sl Via Saldini, 50, I-20133 Milano, Italy}\\
{\rm E-mail:~~\tt rocca@mat.unimi.it}\\}
\maketitle
\date{}
\vspace{-.4cm}

\noindent {\bf Abstract.} In this paper we present a new
thermodynamically consistent phase transition model describing the
evolution of a liquid substance, e.g., water, in a rigid container
$\Omega$ when we  freeze the container. Since the density
$\varrho_{2}$ of ice with volume fraction $\beta_{2}$, is lower
than the density $\varrho_{1}$ of water with volume fraction
$\beta_{1}$, experiments - for instance the freezing of a glass
bottle filled with water - show that the water pressure increases
up to the rupture of the bottle. When the container is not
impermeable, freezing may produce a non-homogeneous material, for
instance water ice or sorbet. Here we describe a general class of
phase transition processes including this example as particular
case. Moreover, we study the resulting nonlinear and singular PDE
system from the analytical viewpoint recovering existence of a
global (in time) weak solution and also uniqueness for some
particular choices of the nonlinear functions involved.

\vspace{.4cm}


\vspace{4mm}

\noindent
{\bf  AMS (MOS) subject clas\-si\-fi\-ca\-tion:}\thspace 80A22, 34A34, 74G25.

\pagestyle{empty}

%\newpage

%%%%%%%%%%%%%%%%%%%%%%%%%%%%%%%%%%%%%%%%%%%%%%%%%%%%%%%%%%%%%%%%%%%
%%%%%%%%%%%%%%%%%%%%%% Caption da ambo i lati %%%%%%%%%%%%%%%%%%%%%

\pagestyle{myheadings}
\newcommand\testopari{\sc Fr\'emond--Rocca}
\newcommand\testodispari{\sc }
\markboth{\testodispari}{\testopari}

%%%%%%%%%%%%%%%%%%%%%%%%%%%%%%%%%%%%%%%%%%%%%%%%%%%%%%%%%%%%%%%%%%%%%%%

\section{Introduction}
\label{intro}

In this paper we consider a phase change process occurring in a
bounded container $\Omega$ of $\RR^3$. More in particular, we
consider a liquid substance (with density $\varrho_1$) in a rigid
and impermeable container $\Omega$ and freeze  it. We
assume that the solid phase has zero velocity, while the liquid
phase velocity is  the  vector ${\bf  U}_{1}$, let $\beta_1$ be the
local proportion of the liquid phase, $\beta_2$ the local
proportion of the solid phase (with density
$\varrho_2<\varrho_1$), and $\teta$ the absolute temperature of
the system. In this setting (cf.~the following Section~\ref{model}
for further explanations on the model), it is possible to  observe macroscopically the formation of
voids, thus we do not have  the condition $\beta_1+\beta_2=1$, but we
only know that
 $\beta_{1}+\beta_{2}\in [0,1]$.  In the following Section~\ref{model} we
introduce the model which will lead to  the following PDE system
we  study - coupling it with suitable initial and boundary
conditions - from the analytical viewpoint in the next
Sections~\ref{exi}--\ref{exiproof}. We search for solutions
$(\teta,\beta_1,\beta_2,p)$ ($p$ standing for the pressure) of the
 Cauchy problem associated with following
PDE system,  coupled with suitable boundary conditions, in the cylindrical domain $Q_T:=\Omega\times(0,T)$,
being $T$ the final time of the process:
\begin{align}\label{eq1}
&k\dt p+\varrho_{1}\dt\beta_{1}+\varrho_{2}\dt\beta_{2}-\Delta p=0,\\
\label{eq2} &\dt(\log\teta)+\frac{\ell}{\teta_c}\dt\beta_{1}
-\Delta(\log\teta)=R,\\
\label{eq3}
&\displaystyle \left(
\begin{array}
[c]{c}%
\displaystyle\dt\beta_{1}\\
\displaystyle\dt\beta_{2}
\end{array}
\right)  -\nu\left(
\begin{array}
[c]{c}%
\Delta\beta_{1}\\
\Delta\beta_{2}%
\end{array}
\right)  -p\left(
\begin{array}
[c]{c}%
\varrho_{1}\\
\varrho_{2}%
\end{array}
\right)  +\left(
\begin{array}
[c]{c}%
\displaystyle-\frac{\ell}{\teta_c}(\teta-\teta_c)\\
0
\end{array}
\right)  +\partial\fhi(\beta_{1},\beta_{2})\ni0.
\end{align}
Here  $\Delta$ represents the Laplacian with respect to the space variables and
we have assumed that $k$, $\ell$ and $\teta_c$ are three
constants denoting, the compressibility of the substance, the
latent heat and the critical temperature of the phase transition
process, respectively, while $R$ stands for an entropy source.
Moreover, $\nu$ represents a non-negative interfacial energy
coefficient and $\varphi$ is a general proper, convex and
lower-semicontinuous function, taking possibly value $+\infty$ in
some points. Hence, $\partial\fhi$ represents here the
subdifferential of Convex Analysis (cf., e.g., \cite{moreau}) and
it takes possibly into account of the restrictions on the phase
field variables $\beta_1$ and $\beta_2$. Indeed, it may be the
indicator function of the set of admissible values of $\beta_1$
and $\beta_2$. If $K=\{(\beta_1,\beta_2)\in
[0,1]\,|\,\beta_1+\beta_2\in [0,1]\}$, then a possible choice in
our analysis is $\fhi(\beta_1,\beta_2)=I_K(\beta_1,\beta_2)$,
which is defined to be $0$ if $(\beta_1,\beta_2)\in K$ and
$+\infty$ otherwise. This is the reason why we have an inclusion
symbol in \eqref{eq3}, since the $\partial\fhi$ may possibly be a
multivalued operator.

We refer to Section~\ref{model} for the detailed description of
of our model, but we prefer to anticipate few words
here on the derivation of equations (\ref{eq1}--\ref{eq3}).
Equation \eqref{eq1} follows from the mass balance  (cf.~also
\cite{fr}) and the equation of macroscopic motion, which is here
assumed to be in agreement with the standard Darcy law
(cf.~\cite{darcy}). Equation \eqref{eq2} represents the rescaled
internal energy balance, which is a consequence of a new approach
developed in \cite{bcf} (and studied also in
\cite{bonetti}--\cite{bcfg2}, \cite{bf}--\cite{br}, and
\cite{fr}). Note that the presence of the logarithmic contribution
in \eqref{eq2} introduce new difficulties in the mathematical
analysis of the problem, due to its singular character, but, once
one has proved that a solution exists in some proper sense, this
immediately entails the positivity of the temperature, which could
otherwise represent a difficult step in the proof (cf., e.g.,
\cite{cspos}). Inclusion \eqref{eq3} stands for  the microscopic
motion equation (cf.~\cite{fremond} and also \cite{bs} and
\cite{visintin} for different approaches) for a proper choice of
the free-energy functional and of the pseudo-potential of
dissipation which, following the theory of \cite{moreau}, accounts
for dissipation in the model (cf.~also the next
Section~\ref{model} for further details on this topic).

In all references  quoted above (but in \cite{fr}) only
phase transition processes in which no voids nor overlapping between phases can occur are studied.
In  that case the proofs mainly rely on the fact that  $\beta_1+\beta_2=1$,
hence the mass balance equation  as well as the effects of the pressure can be
neglected and so the system reduces to the coupling of equation \eqref{eq2} (with only one phase variable, say
$\chi=\beta_1-\beta_2$) and \eqref{eq3}.
More in particular, in \cite{bcf} the model coupling these two equations and accounting also
for some memory effects in \eqref{eq2} is studied
and an existence - of weak solution - result is proved in case $\fhi=I_{[0,1]}$, while
uniqueness remains still an open problem. In \cite{bf} and \cite{bcfg1,bcfg2}
uniqueness of solution is   established and also some investigation on
the long-time behaviour of solutions is performed respectively in case $\nu=0$ and without any thermal memory
for a system in which  a different heat flux law in \eqref{eq2} - leading to an internal energy balance
with $\teta$ instead of $\log\teta$ inside the Laplacian - is taken into account.
Moreover, in \cite{br} (and in \cite{bonetti} for the case $\nu=0$ and
without any thermal memory contribution)
the long-time behaviour of the solutions-trajectories to an extension of the model
presented in \cite{bcf} is studied, while in \cite{bfr} a more general approach to
this kind of processes - leading to more general type of nonlinearities in
the equations - is performed. In particular in \cite{bfr}
the $\log\teta$-term in \eqref{eq2} is substituted by a more general
nonlinear function (maybe also multivalued) including the logarithmic function as a
particular case. An existence-long-time behaviour of solution result is established
in \cite{bfr}. Only in \cite{fr} - instead -  a first  attempt to approach phase-field systems including the possibility of having voids
is developed. In \cite{fr} indeed we study a system coupling the phase equation \eqref{eq3} (including also
a viscosity term), the temperature equation
\eqref{eq2}, a general mass balance equation, and a quasi-static standard stress-strain relation, including viscosity
effects but no accelerations. We establish in \cite{fr}  a global well-posedness result  for the resulting
PDE system coupled with suitable initial and boundary conditions.

Here, under the assumption that
the solid phase velocity is 0 and by means of the Darcy law (cf.~\cite{darcy} and also \eqref{sorbet6} below),
the mass balance equation can be rewritten as \eqref{p1}. The main mathematical
difficulties  stand in the
nonlinear and singular coupling between \eqref{p2} and \eqref{p3}. Our approach to this
system is to regularize the $\log$-term  in \eqref{eq2} by introducing its Lipschitz-continuous approximation. Hence, we
first solve the regularized problem in Section~\ref{proofpreli}, getting also uniqueness of solutions,
then, in the last Section~\ref{exiproof} we prove our main result (stated in Section~\ref{exi}) entailing
global (in time) existence of (weak) solutions to (\ref{eq1}--\ref{eq3}) coupled with suitable initial
and boundary conditions. The problem of uniqueness for system (\ref{eq1}--\ref{eq3})
is still open (cf.~also Remark~\ref{unioss} below)  and this is mainly due to the lack of regularity of the $\teta$-component of the solution $(\teta,\beta_1,\beta_2,p)$ to (\ref{eq1}--\ref{eq3}).

%%%%%%%%%%%%%%%%%%%%%%%%%%%%%%%%%%%%%%%%%%%%%%%%%%%%%%%%%%%%%%%%%%%%%%%%%%%%%%%%%

\section{The model}
\label{model}

In this Section we detail the derivation of the model leading to
the PDE system (\ref{eq1}--\ref{eq3}) already introduced in
Section~\ref{intro}. Hence, let us consider a liquid substance,
e.g., water, in a rigid and impermeable container $\Omega$ and
freeze the container. Since the density $\varrho_{2}$ of ice with
volume fraction $\beta_{2}$ is lower than the density $\varrho_{1}$ of
water with volume fraction $\beta_{1}$, it seems impossible that
water transforms into ice (or, more in general, into a mixture of
ice and water): the water remains liquid even at low temperature.
Experiments - for instance the freezing of a glass bottle filled
with water - show that the water pressure increases up to the
rupture of the bottle. When the container is not impermeable,
freezing may produce an non homogeneous material, for instance
water ice or sorbet. In this setting, it is possible to observe
macroscopically the presence of
voids and thus we do not have the condition $\beta_1+\beta_2=1$, but we
can only deduce that
 $0\leq\beta_{1}+\beta_{2}\leq1$.

We assume the ice has zero velocity, while the water velocity is
represented by the vector ${\bf  U}_{1}$. For the sake of simplicity, we assume the
small perturbation assumption (cf.~\cite{germain}). In this
setting the mass balance is
\begin{align}\no
&\beta_{1}\dt\varrho_{1}+\beta_{2}\dt\varrho_{2}+\varrho_{1}\dt\beta_{1}
+\varrho_{2}\dt\beta_{2}+\varrho_{1}\beta_{1}\dive{\bf U}_1=0,\\
\label{sorbet0} &-\varrho_{1}\beta_{1}{\bf U}_1\cdot{\bf
n}=\varpi,
\end{align}
where $\varpi$ is the liquid water intake and by ${ \bf n}$ we
denote the outward unitary normal vector to $\partial\Omega$. We
have assumed the substance to be compressible, even in some cases
this compressibility is very low. Moreover, the equation of motion
for the macroscopic motion of the water can be written as
\begin{align}\no
&\varrho_{1}\dt\left(  \beta_{1}{\bf U}_{1}\right) =\operatorname{div}%
\sigma_{1}+{\bf f}_{1},\\
\label{sorbet1}
&\sigma_{1}{\bf n}={\bf g}_{1},
\end{align}
where $\sigma_{1}$ is the water stress, ${\bf f}_{1}$ is the
action of the exterior on the liquid water, for instance the
friction of the water on the solid ice phase, and ${\bf g}_{1}$ is
the boundary exterior action on the liquid water.

For the microscopic  forces  responsible for the ice water phase
change, the equation of motion is (cf.~\cite{fremond})
\begin{equation}\no
-\left(
\begin{array}
[c]{c}%
B_{1}\\
B_{2}%
\end{array}
\right)  +\operatorname{div}\left(
\begin{array}
[c]{c}%
{\bf H}_{1}\\
{\bf H}_{2}%
\end{array}
\right)  =0.
\end{equation}
This vector equation is the equation of microscopic motion and it
deals with the works produced by the microscopic motions occurring in the phase changes.
The terms $B_1$, $B_2$, ${\bf H}_1$, and ${\bf H}_2$ are new internal forces
responsible for phase transitions phenomena. In particular $B_1$
 represents the energy density per unit of $\beta_1$, $B_2$ represents
 the energy density per unit of $\beta_2$, while $H_1$ and $H_2$ represent densities of energy
 flux. If we assume that no exterior source of work is present, we also have
\begin{equation}
\label{sorbet2}%
\left(
\begin{array}
[c]{c}%
{\bf H}_{1}\\
{\bf H}_{2}%
\end{array}
\right)  \cdot{\bf n}=0.
\end{equation}
The internal energy balance
reformulated in terms of the entropy $s$ is (cf.~\cite{bcf} for
its derivation)
\begin{align}\label{sorbetent}
\dt s+\operatorname{div}{\bf Q}=&\,R,
\end{align}
where we have denoted by ${\bf Q}$ the entropy flux, by $R$ the
entropy source, $r=R\teta$ being the heat source, and by $\betab$
the column vector of components $\beta_1,\,\beta_2$. No high order
dissipative contribution is involved on the right hand side due to
the small perturbation assumption.

The internal forces are split into dissipative forces indexed by
$^{d}$
and non dissipative interior forces indexed by $^{nd}$.
The non dissipative forces are defined through the water free energy  $\Psi$
which takes into account the internal constraint on the state
quantities, for instance%
\begin{equation}
(\beta_{1},\beta_{2})\in K=\{(\gamma_1,\,\gamma_2)\in
[0,1]\,|\,\gamma_{1}+\gamma_{2}\in[0,1] \}.\nonumber
\end{equation}
We choose the following form for the free energy:
\begin{align}\no
\Psi(\teta,\varrho_{1},\varrho_{2},\beta_{1},\beta_{2},\nabla\beta_{1},\nabla\beta_{2})=&\,-c_0\teta(\log\teta
-1)-\frac{\beta_{1}\ell}{\teta_c}(\teta-\teta_c)+I_{K}(\beta_{1},\beta_{2})\\
\label{psi} &+\frac{\nu}{2}\left(
 |\nabla\beta_{1}|^2+
 |\nabla\beta_{2}|^2\right)+\beta_{1}f_{1}(\varrho_{1})+\beta_{2}f_{2}(\varrho_{2}),
\end{align}
where $\teta_{c}$ is the critical temperature of the system,
$\ell$ denotes the latent heat of the melting solidification
process, $c_0$ the positive specific heat, $\nu$ a
positive interfacial energy coefficient, and $I_K$ stands for the
indicator function of the set $K$, it takes value $0$ on $K$ and
$+\infty$ otherwise.  Because we have assumed the substance to be
compressible, the free energy depends on the densities. Indeed,
$f_1$ and $f_2$  account, respectively, for the dependence on the
densities of the free energy related to the liquid and the solid
phases. Thus we have
\begin{align}
&\displaystyle s=-\frac{\partial\Psi}{\partial\teta}=c_0\log\teta
+\frac{\ell\beta_1}{\teta_c},\label{sorbetropi}
\end{align}
and we define the non dissipative internal forces as:
\begin{align}
&\displaystyle\left(
\begin{array}
[c]{c}%
B_{1}^{nd}\\
B_{2}^{nd}%
\end{array}
\right) = \left( \begin{array}
[c]{c}%
\displaystyle\frac{\partial\Psi}{\partial\beta_{1}}\\
\displaystyle\frac{\partial\Psi}{\partial\beta_{2}}%
\end{array}
\right)  \in \left(
\begin{array}
[c]{c}%
\displaystyle-\frac{\ell}{\teta_c}(\teta-\teta_c)\\
0
\end{array}
\right)  +\partial I_{K}(\beta_{1},\beta_{2}),\label{sorbet4nd}\\
&\left(
\begin{array}
[c]{c}%
{\bf H}_{1}^{nd}\\
{\bf H}_{2}^{nd}%
\end{array}
\right) = \left( \begin{array}
[c]{c}%
\displaystyle\frac{\partial\Psi}{\partial (\nabla\beta_{1})}\\
\displaystyle\frac{\partial\Psi}{\partial (\nabla\beta_{2})}%
\end{array}
\right)
 =\nu\left(
\begin{array}
[c]{c}%
 \nabla\beta_{1}\\
 \nabla\beta_{2}%
\end{array}
\right)  ,\label{sorbet5nd}\\
&\displaystyle\left(
\begin{array}
[c]{c}%
P_{1}^{nd}\\
P_{2}^{nd}%
\end{array}
\right) = \left( \begin{array}
[c]{c}%
\displaystyle\frac{\partial\Psi}{\partial\varrho_{1}}\\
\displaystyle\frac{\partial\Psi}{\partial\varrho_{2}}%
\end{array}
\right) =\left(
\begin{array}
[c]{c}%
 \beta_{1}\displaystyle\frac{\partial f_{1}}{\partial\varrho_{1}}\\
 \beta_{2}\displaystyle\frac{\partial f_{2}}{\partial\varrho_{2}}%
\end{array}
\right) , \label{sorbet6nd}
\end{align}
with an abuse of notation for the derivatives
$\partial\Psi/\partial \beta$ and where $\partial I_K$ represents
the subdifferential of $I_K$ in the sense of Convex Analysis (cf.,
e.g., \cite{brezis} or \cite{moreau}). Because the free energy
does not depend on deformations, we have
\begin{equation}
\sigma_{1}^{nd}=0.\label{sorbetmf1}
\end{equation}
The dissipative forces are defined with
pseudopotential of dissipation $\Phi$ which takes into account the
internal constraints related to the velocities and ensures that
the second law of thermodynamics is satisfied. The mass balance represents
such a constraint (cf.~also \cite{fr} for further comments on this
choice).
For instance, we may choose%
\begin{align}\no
\Phi\Big(\dt\varrho_{1},\dt\varrho_{2},\dt\beta_{1},\dt\beta_{2}
,& \nabla(\dt\beta_{1}), \nabla(\dt\beta_{2}),D({ \bf U }_{1}),
\nabla\teta\Big)\\
\no
=&\,I_{0}\left(\beta_{1}\dt\varrho_{1}+\beta_{2}\dt\varrho_{2}+\varrho_{1}\dt\beta_{1}+\varrho_{2}\dt\beta_{2}
+\varrho_{1}\beta_{1}\dive{ \bf U }_1\right)\\
\no
&+\frac{\mu}{2}\left( \left|\dt\beta_{1}\right|^2+
\left|\dt\beta_{2}\right|^2\right)
+\frac{\lambda}{2\teta}|\nabla\teta|^{2},
\end{align}
where $I_{0}$ denotes the indicator function of the origin of
$\RR$, $\mu$ the phase change viscosity and $\lambda>0$ denotes
the heat conductivity of the system, while the operator $D$
represents the linearized symmetric strain tensor $D _{ij}(
\mathbf{U}):=(U_{i,j}+U_{j,i})/2$, where $i,j=1,2,3$ and with the
commas we mean space derivatives. The dissipative forces can be written as:%
\begin{align}\no
(P_{1}^{d},\,P_{2}^{d},\,B_{1}^{d},\,B_{2}^{d},\,{\bf H}_{1}^{d},&\,{\bf
H}_{2}^{d},\,\sigma_{1}^{d},\,{\bf Q}^d)\\
\no
&\in
\partial\Phi\left(\dt\varrho_{1},\,\dt\varrho_{2},\,\dt\beta_{1},\,\dt\beta_{2},\, \nabla(\dt\beta_{1}),
\,\nabla(\dt\beta_{2}),\,D({ \bf U }_{1}),\, \nabla\teta\right),
\end{align}
where $\partial \Phi$ represents the subdifferential of $\Phi$ and it
results
\begin{align}
&{\bf Q}={\bf Q}^d=-\frac{\partial\Phi}{\partial\nabla\teta}=-\lambda\nabla(\log\teta),\label{sorbet3}\\
&{\bf H}_{1}^{d}
=\frac{\partial\Phi}{\partial(\nabla\dt\beta_1)}=0,
\ {\bf H}_{2}%
^{d}=\frac{\partial\Phi}{\partial(\nabla\dt\beta_2)}=0,\label{sorbet4}\\
&\displaystyle\left(
\begin{array}
[c]{c}%
P_{1}^{d}\\
P_{2}^{d}\\
B_{1}^{d}\\
B_{2}^{d}\\
\sigma_{1}^{d}%
\end{array}
\right) = \left(
\begin{array}
[c]{c}%
\displaystyle\frac{\partial\Phi}{\partial(\dt\varrho_1)}\\
\displaystyle\frac{\partial\Phi}{\partial(\dt\varrho_2)}\\
\displaystyle\frac{\partial\Phi}{\partial(\dt\beta_1)}\\
\displaystyle\frac{\partial\Phi}{\partial(\dt\beta_2)}\\
\displaystyle\frac{\partial\Phi}{\partial D({\bf
U}_1)}%
\end{array}
\right) \in \left(
\begin{array}
[c]{c}%
0\\
0\\
\mu\dt\beta_{1}\\
\mu\dt\beta_{2}\\
0%
\end{array}
\right) +\partial I_{0}(0)\left(
\begin{array}
[c]{c}%
\beta_{1}\\
\beta_{2}\\
\varrho_{1}\\
\varrho_{2}\\
\varrho_{1}\beta_1\mathbf{1}%
\end{array}
\right),\nonumber
\end{align}
with an abuse of notation for the derivatives
$\partial\Phi/\partial(\dt X)$ and where $\mathbf{1}$ stands for the unit
tensor. By denoting $-p$ an element of $\partial I_0$, the last
relations can be written as:
\begin{align}
&P_{1}^{d}=-p\beta_{1},\ P_{2}^{d}=-p\beta_{2},\label{sorbet6bis}\\
&B_{1}^{d}=\mu\dt\beta_{1}-p\varrho_{1},\ B_{2}^{d}=\mu\dt\beta_{2}-p\varrho_{2},\label{sorbet4bis}\\
&\sigma_{1}^{d}=-p\varrho_{1}\beta_{1}\mathbf{1}.\label{sorbet3bis}%
\end{align}
The quantity $p\varrho_{1}\beta_{1}$ turns out to be the pressure in
the liquid phase and using constitutive law (\ref{sorbet3bis}) in
the equation for the macroscopic motion (\ref{sorbet1}) we get
\[
\varrho_{1}\dt\left(  \beta_{1}{ \bf U }_{1}\right)=- \nabla%
\left(p\varrho_{1}\beta_{1}\right)+{\bf f}_{1}.
\]
We may assume that the acceleration is negligible and that the volume
exterior force ${\bf f}_{1}$ results mainly from friction on the
solid
phase and that this force is proportional to the relative velocity through the relation
\[
{\bf f}_{1}=-\frac{\left(  \varrho_{1}\beta_{1}\right)  ^{2}}{m}{ \bf U }_{1}.
\]
This is mainly due to the fact that the force has to be ${\bf 0}$ when there is no water.

\beos
The density of force ${\bf f}_{1}$ with respect to the actual volume
of water $\beta_{1}d\Omega$ is
\[
\frac{{\bf f}_{1}}{\beta_{1}}=-\frac{\varrho_{1} ^{2}\beta_{1}}{m}{
\bf U }_{1}.
\]
It is proportional to $\beta_{1}$ and to the liquid water momentum
$\varrho_{1}\beta_{1}{\bf U }_{1}$ as one may expect.
 \eddos

\noindent
With this choice, from \eqref{sorbet1} we get%
\[
\frac{\left(  \varrho_{1}\beta_{1}\right)  ^{2}}{m}{ \bf U }_{1}%
=-\nabla \left(p\varrho_{1}\beta_{1}\right),
\]
which (within the small perturbation theory (cf.~\cite{germain})) can be rewritten as:
\begin{equation}
\varrho_{1}\beta_{1}{ \bf U }_{1}=-m\nabla p,\label{sorbet6}%
\end{equation}
where the parameter $m$ denotes the mobility of the water in agreement
with the Darcy law (cf. \cite{darcy}).

Since the densities $\varrho_1$ and $\varrho_2$ are not related to
a particular internal force, their constitutive laws are
(cf.~\cite[(3.25), p.~11]{fremond})
\begin{equation}\no
P_{1}^{nd}+P_{1}^{d}=0, P_{2}^{nd}+P_{2}^{d}=0,
\end{equation}
which, with (\ref{sorbet6nd}) and \eqref{sorbet6bis}, can be written as:
\begin{equation}
\displaystyle\frac{\partial f_{1}}{\partial\varrho_{1}}-p=0,
\displaystyle\frac{\partial f_{2}}{\partial\varrho_{2}}-p=0.
\end{equation}
Within the small perturbation assumption, we have
\begin{equation}\label{eqprho}
p=k_{1}(\varrho_{1}-\varrho_{1}^{r}),
p=k_{2}(\varrho_{2}-\varrho_{2}^{r}),
\end{equation}
where the $k_{i}$ ($i=1,\,2$) are the compressibilities of the liquid and
solid phases, respectively, and the $\varrho_{i}^{r}$ ($i=1,\,2$) are reference densities.
 Let us note that the densities,
which depend on the pressure $p$ through \eqref{eqprho}, may also depend on
the temperature via a thermal expansion coefficient, giving
$p=k_{i}(\varrho_{i}-\varrho_{i}^{r})+k_{i}^{th}(\teta-\teta_{c})$ in
place of \eqref{eqprho}. However, for the time being, we do not introduce this dependence.

The predictive theory equations within the small perturbation
assumption result from  the mass balance (\ref{sorbet0})  with \eqref{eqprho}, the equation of
motion (\ref{sorbet1}--\ref{sorbet2}),  the entropy balance
(\ref{sorbetent}) and  the  constitutive laws
(\ref{sorbetropi}--\ref{sorbet4bis})  with   relation (\ref{sorbet6}).  Hence, they
can be written as follows
\begin{align}\label{sorbet7}
&k\dt p+\varrho_{1}\dt\beta_{1}+\varrho_{2}\dt\beta_{2}-m\Delta p=0,\\
\label{sorbet7p}
&c_0\dt(\log\teta)+\frac{\ell}{\teta_c}\dt\beta_{1}
-\lambda\Delta(\log\teta)=R,\\
\label{sorbet7s}
&\displaystyle \mu\left(
\begin{array}
[c]{c}%
\displaystyle\dt\beta_{1}\\
\displaystyle\dt\beta_{2}
\end{array}
\right)  -\nu\left(
\begin{array}
[c]{c}%
\Delta\beta_{1}\\
\Delta\beta_{2}%
\end{array}
\right)  -p\left(
\begin{array}
[c]{c}%
\varrho_{1}\\
\varrho_{2}%
\end{array}
\right)  +\left(
\begin{array}
[c]{c}%
-\frac{\ell}{\teta_c}(\teta-\teta_c)\\
0
\end{array}
\right)  +\partial I_{K}(\beta_{1},\beta_{2})\ni0,
\end{align}
where $\displaystyle k=\left(\frac{\beta_{1}^{0}}{k_{1}}+\frac{\beta_{2}^{0}}{k_{2}}\right)$
is the compressibility of the mixture,
being $\beta_1^0+\beta_2^0$ the reference value of the material volume fraction.
This set of equation is completed by the following boundary conditions on $\partial\Omega$%
\begin{align}\label{bobeta}
&\frac{\partial\beta_{1}}{\partial {\bf n}}=0,\ \frac{\partial\beta_{1}}{\partial
{\bf n}}=0,\\
\label{boteta}
&\lambda\frac{\partial\teta}{\partial {\bf n}}+\alpha_{\teta}(\teta-\teta
_{ext})=0,
\end{align}
where $\alpha_{\teta}$ is the thermal conductivity of the container boundary.
Concerning the hydraulic boundary conditions, we may choose%
\begin{equation}\label{bop}
m\frac{\partial p}{\partial {\bf n}}+\alpha_{p}(p-p_{ext})=0,
\end{equation}
where $\alpha_{p}=0$ if the boundary is watertight, $\alpha_{p}=\infty$
if the boundary is connected to a water supply with pressure
$p_{ext}$ (in this case the boundary condition is $p=p_{ext}$).
This condition means that the flow is proportional to the
difference of pressure $p-p_{ext}$, in case the container is permeable.
In case, instead, of a metallic recipient, the condition accounts for the
deformation of the metal due to the pressure.

\beos Due to the small perturbation assumption, densities
$\varrho_{1}$ and $\varrho_{2}$ in equations (\ref{sorbet7}) and
(\ref{sorbet7s}) are constant.
 \eddos
 \beos\label{condlim} More sophisticated boundary condition
result from more sophisticated physical boundary properties. The
boundary of the container can be semi-permeable, in this case we
have
  \[
m\frac{\partial p}{\partial {\bf
n}}-\alpha_{p}(p-p_{ext})^{-} = 0,
\]
where $(p-p_{ext})^{-}=\sup(0,-(p-p_{ext}))$ is the negative part
of $p-p_{ext}$.
This boundary condition means that when the pressure is lower than
the exterior pressure $p_{ext}$, water flows inside the container
but when the pressure is larger than the exterior pressure, no
water flows outside.

The pressure may also be controlled on the
boundary by the following relation
\[
m\frac{\partial p}{\partial {\bf n}}+\partial
I_{-}(p-p_{ext})\ni0.
\]
Water flows outside the container in order to maintain the
pressure lower than the outside pressure $p_{ext}$ which may be
the atmospheric pressure. \eddos

The initial conditions in $\Omega$ are%
\[
\beta_{1}(x,0)=\beta_{1}^{0}(x),\ \beta_{2}(x,0)=\beta_{2}^{0}(x),\ \teta
(x,0)=\teta_0(x),\ p(x,0)=p_0(x) .
\]
This set of partial differential equations allow to compute pressure $p(x,t)$,
liquid water and ice volume fractions $\beta_{1}(x,t)$, $\beta_{2}(x,t)$ and
temperature $\teta(x,t)$ depending on the external actions resulting from the
exterior pressure $p_{ext}(x,t)$ and exterior temperature $\teta_{ext}(x,t)$
and rate of heat production $\teta_cR(x,t)$ very often equal to $0 $. In an
engineering situation indeed the governing action is the exterior temperature cooling
and heating the container.

\subsection{Some examples}
\label{examples}

\paragraph{1. Freezing water in an impermeable container.}

Suppose to  have $\beta_{1}^{0}=1,\ \beta_{2}^{0}=0,\ \teta_0>0$ and that the
boundary condition for the pressure is (\ref{bop}) with
$\alpha_{p}=0$, i.e.:
\[
\frac{\partial p}{\partial {\bf n}}=0.
\]
Hence, the container is completely cooled. In order to look for closed form
solutions we assume all the quantities are homogeneous. Assuming
the temperature to be known, we focus on equation (\ref{sorbet7s})
which describes the phase change. Because of homogeneity, it is%
\[\displaystyle
\mu \left(
\begin{array}
[c]{c}%
\dt\beta_{1}\\
\dt\beta_{2}
\end{array}
\right)  -p\left(
\begin{array}
[c]{c}%
\varrho_{1}\\
\varrho_{2}%
\end{array}
\right)  +\left(
\begin{array}
[c]{c}%
-\frac{\ell}{\teta_c}(\teta-\teta_c)\\
0
\end{array}
\right)  +\partial I_{K}(\beta_{1},\beta_{2})\ni0.
\]
When the Celsius temperature $\teta-\teta_c$ is positive, the
solution of this equation is $\beta _{1}=1$ and $\beta_{2}=0$.
When the Celsius temperature becomes negative, the phase change
which is expected does not occur. Indeed due to the mass balance,
to
values $\beta_{1}=1$ and $\beta_{2}=0$ and to condition $\beta_{1}+\beta_{2}%
\leq1$, we should have, neglecting the effect of the
compressibility, the following relations:
\[
\varrho_{2}\dt\beta_{2}+\varrho_{1}\dt\beta_{1}=0,\dt\beta_{1}<0,\dt\beta
_{2}+\dt\beta_{1}\leq0.
\]
Note that the two first conditions imply%
\[
\dt\beta_{2}+\dt\beta_{1}=\left(1-\frac{\varrho_{1}}{\varrho_{2}}\right)\dt\beta_{1}>0,
\]
which contradicts the third condition. Let us check that the water
remains liquid ($\dt\beta_{1}=0$) even if Celsius temperature is negative. We must have%
\[
-p\left(
\begin{array}
[c]{c}%
\varrho_{1}\\
\varrho_{2}%
\end{array}
\right)  +\left(
\begin{array}
[c]{c}%
-\frac{\ell}{\teta_c}(\teta-\teta_c)\\
0
\end{array}
\right)  +\partial I_{K}(1,0)\ni0,
\]
with%
\[
\partial I_{K}(1,0)=\left(
\begin{array}
[c]{c}%
P\\
Q
\end{array}
\right)  ,\ P\geq0,\ Q\leq P,
\]
or%
\begin{gather*}
p\varrho_{1}+\frac{\ell}{\teta_c}(\teta-\teta_c)\geq0,\\
0\leq
p(\varrho_{1}-\varrho_{2})+\frac{\ell}{\teta_c}(\teta-\teta_c).
\end{gather*}
Since $\varrho_{1}-\varrho_{2}>0$, these conditions may be satisfied. In this case, the minimum
value of the pressure is%
\[
p=-\frac{\ell}{(\varrho_{1}-\varrho_{2})\teta_c}(\teta-\teta_c).
\]
Thus, when cooling the container, water remains liquid and pressure
increases. This is in agreement with experiment: a glass container
filled of water explodes due to the pressure increase when
freezing.

\paragraph{2. Freezing of a water emulsion in an impermeable
container.}

Suppose to have $0<\beta_{1}^{0}<1,\ \beta_{2}^{0}=0,\ \teta_0>0$, hence, the emulsion (a
mixture of voids and water) is cooled. Assuming an homogeneous evolution, it
may be seen that the water freezes till either $\beta_{1}=0$, (in this case
$\beta_{1}^{0}<\varrho_{2}/\varrho_{1}$ and the cooling results in a mixture of ice
and voids), or $\beta_{1}+\beta_{2}=1$, (in this case $\beta_{1}^{0}\geq
\varrho_{2}/\varrho_{1}$ and the cooling results in a mixture of ice and liquid
water, $\beta_{1}=(\varrho_{1}\beta_{1}^{0}-\varrho_{2})/(\varrho_{1}-\varrho_{2})$,
$\beta_{2}=(\varrho_{1}-\varrho_{1}\beta_{1}^{0})/(\varrho_{1}-\varrho_{2})$). This
mixture is a water ice or sorbet.

\paragraph{3. The ice in a glacier.}

Assume to have an homogeneous mixture of ice and water without voids at an
equilibrium,
\[
\beta_{1}+\beta_{2}=1,\ \dt\beta_{1}= \dt\beta_{2}=0,
\]
then, equations (\ref{sorbet7}--\ref{sorbet7s}) give
\[
-p\left(
\begin{array}
[c]{c}%
\varrho_{1}\\
\varrho_{2}%
\end{array}
\right)  +\left(
\begin{array}
[c]{c}%
-\frac{\ell}{\teta_c}(\teta-\teta_c)\\
0
\end{array}
\right)  +\left(
\begin{array}
[c]{c}%
P\\
P
\end{array}
\right)  =0,
\]
with $P\geq0$. Hence, we get
\[
-p(\varrho_{1}-\varrho_{2})-\frac{\ell}{\teta_c}(\teta-\teta_c)=0,\
p\geq0.
\]
In this situation the phase change temperature depends on the pressure. It is
known that this is the case for the ice-water phase change. Because $\varrho
_{1}>\varrho_{2}$ the phase change temperature decreases when the pressure
increases. A consequence of this property is that the ice melts at the bottom
of a glacier and lubricates the rock ice contact surface allowing the downhill
motion of the glacier.
%%%%%%%%%%%%%%%%%%%%%%%%%%%%%%%%%%%%%%%%%%%%%%%%%%%%%%%%%%%%%%%%%%%%%%%%

%%%%%%%%%%%%%%%%%%%%%%%%%%%%%%%%%%%%%%%%%%%%%%%%%%%%%%%%%%%%%%%%%%%
%%%%%%%%%%%%%%%%%%%%%%%%%%%% exiuni  %%%%%%%%%%%%%%%%%%%%%%%%%%%%%%%%

\section{Variational formulation and main results}
\label{exi}

In the following part of the paper we want to deal with a suitable generalization
of the PDE system (\ref{sorbet7}--\ref{sorbet7s}) introduced in the previous
Section~\ref{model}. Hence, we first detail in this section the assumptions
on the data and on the general proper, convex, and  \lsc{}
mapping $\varphi: \RR^2 \to[0,+\infty]$, which will generalize the
role played by the indicator function $I_K$ in the inclusion \eqref{sorbet7s}
introduced in the previous Section~\ref{model}, and then we state
a suitable variational formulation of
the generalized PDE system obtained and our main results.

In order to do that, we first introduce the Hilbert
triplet $(V,H,V')$ where
\begin{equation}\label{spazi}
H:=L^2(\Omega)\quad\hbox{and } V:=H^1(\Omega),
\end{equation}
and $\Omega\subset\RR^3$ is a bounded and  connected domain with Lipschitz continuous
boundary $\Gamma:=\partial\Omega$. Let $T$ be a positive final time of the process and denote by
$Q_t:=\Omega\times (0,t)$ and $\Sigma_t:=\Gamma\times(0,t)$, $t\in (0,T]$.
Then, we identify, as usual, $H$ (which stands either for the space $L^2(\Omega)$ or for $(L^2(\Omega))^3$ or for
$(L^2(\Omega);\RR^2)$) with its dual space $H'$, so that $V\hookrightarrow H\hookrightarrow V'$ with dense
and continuous embeddings. Moreover, we denote by $\Vert\cdot\Vert_X$ the norm in some space $X$, by $(\cdot,\cdot)$
the scalar product in $H$,  and by
$\duav{\cdot,\cdot}$ the duality pairing between $V$ and $V'$.
Hence, for any $\zeta\in V'$, set
\bega{defiwo}
  \zeta_\Omega:=\frac{1}{|\Omega|}\duav{\zeta,1},\\
 \label{defiV}
  {\cal V}':=\{\zeta\in V': \zeta_\Omega=0\},
   \qquad {\cal V}:=V\cap {\cal V}'.
\end{gather}
The above notation $\cal{V}'$ is just suggested for the
sake of convenience; indeed, we mainly see ${\cal V}$, ${\cal V}'$ as
(closed) subspaces of $V$, $V'$, inheriting their norms,
rather than as a couple of spaces in duality.

>From now on, for simplicity and without any loss of generality, we  suppose that the coefficients in
(\ref{sorbet7}--\ref{sorbet7s})
$k=m=c_0=\lambda=\mu=\ell/\teta_c=1$.
Next, in order to give a precise formulation of our problem (\ref{sorbet7}--\ref{sorbet7s}),
 we define here the realization of the Laplace operator
with Neumann homogeneous boundary conditions, that is the operator
\begin{align}\label{laplneu}
&{\cal B}:V\to V',\quad\duav{{\cal B}u, v}=\io\nabla u\cdot\nabla v \quad u,v\in V.
\end{align}
Clearly, ${\cal B}$ maps $V$ onto ${\cal V}'$
and its restriction to ${\cal V}$ is an isomorphism of ${\cal V}$
onto ${\cal V}'$.

Define $W:=\{w\in H^2(\Omega):\quad\dn w=0\quad\hbox{on }\Gamma\}$, where $\dn$
is the derivative with respect to~the outward normal derivative to $\Gamma$ and
 make the following assumptions on the data.
\behy\label{hyp1} We assume that $\varrho_1>\varrho_2>0$, $\nu\geq 0$,
and suppose that
\begin{align}\label{ipophi}
&\varphi:\,\RR^2\to [0,+\infty]\hbox{ is proper, convex, lower semicontinuous and}\quad\varphi({\bf 0})=0,\\
\label{iponu}
&{\cal D}(\varphi) \hbox{ bounded if }\nu=0,\\
\label{gammaexp}
&\gamma(r):=\exp(r)\quad\hbox{for }r\in \RR,\\
\label{datoteta}
&\teta_0\in L^1(\Omega),\quad\teta_0>0\quad\hbox{a.e. in }\Omega,\quad w_0:=\gamma^{-1}(\teta_0)\in H,\\
\label{datochi}
&\betab_0=(\beta_1^0,\beta_2^0)\in {\cal D}(\varphi),\quad \nu\beta_1^0,\,\nu\beta_2^0\in V,\quad p_0\in V,\\
\label{sorg2}
&R\in L^2(Q_T)\cap L^1(0,T;L^\infty(\Omega)),\quad \Pi\in L^{\infty}(\Sigma_T).
\end{align}
Then, we introduce the functions ${\cal R}\in\LDVp$ such that
\begin{align}\label{sorg3}
&\duav{{\cal R}(t),v}=\io R(t)v+\int_{\partial\Omega} \Pi(t) v_{|\partial\Omega}\quad v\in V,\quad
\hbox{for a.e. } t\in [0,T].
\end{align}
\enhy
Moreover, we denote by $\partial\varphi$ the subdifferential in $\RR^2\times\RR^2$ of the convex analysis
and note that $\partial\varphi$ is maximal monotone
and that $\partial\varphi({\bf 0})\ni{\bf 0}$ (see, e.g., \cite{barbu} and \cite{brezis} for the general theory).
The same symbol $\partial\varphi$ will be
used for the maximal monotone operators induced on $L^2$-spaces.

\beos\label{boucond}
Here in our analysis we have chosen to treat the Neumann homogeneous boundary conditions
on $p$, that is, we consider \eqref{bop} with $\alpha_p=0$. It is not difficult, however
to treat the boundary conditions \eqref{bop} with $\alpha_p>0$
with the same techniques used here. Regarding, instead, the absolute temperature $\teta$,
we work here with Neumann non-homogeneous boundary conditions on $\log\teta$, which means
we consider \eqref{boteta} in case $\alpha_\teta>0$ and $\teta_{ext}=0$. Other type of boundary conditions
like Dirichlet non-homogeneous boundary conditions on $\log\teta$ could also be
taken into account without any difficulties in our analysis. Regarding instead the
more sophisticated boundary conditions on $p$ detailed in Remark~\ref{condlim}
a more careful analysis should be performed in order to include them in our results, but
we do not want to face this problem in this paper.
\eddos

Then, we are ready to introduce the variational
formulation of our problem as follows.

\medskip
\noindent
{\sc Problem (P).} Find $(w,\,\beta_1,\beta_2)$ and $(\teta,\,\xi_1,\xi_2,\,p)$ with the regularities
\begin{align}
\label{reg2}
&w\in \HUVp\cap\LDV,\quad\teta\in L^{5/3}(Q_T)\\
\label{reg3}
&\beta_1,\beta_2\in \HUH,\quad \nu\beta_1,\,\nu\beta_2\in \LIV\cap L^{5/3}(0,T;W^{2,5/3}(\Omega)),\\
\no
&\betab=(\beta_1,\beta_2)\in {\cal D}(\varphi)\quad\hbox{a.e. in }Q_T,\\
\label{reg4}
&\xi_1,\xi_2\in L^{5/3}(Q_T),\quad p\in \CZV\cap\LDW\cap\HUH,
\end{align}
satisfying
\begin{align}\label{p1}
&\dt p+\varrho_1\dt\beta_1+\varrho_2\dt\beta_2+{\cal B}p=0\quad\hbox{ a.e. in }Q_T,\\
\label{p2}
&\dt w
+\dt\beta_1+{\cal B}w={\cal R}\quad\hbox{in }V'\hbox{ a.e. in }[0,T],\\
\label{p3}
&\displaystyle\dt{\betab}
+\nu\begin{pmatrix}{\cal B}\beta_1\\{\cal B}{\beta_2}\end{pmatrix}+\xib-
p\begin{pmatrix}\varrho_1\\\varrho_2\end{pmatrix}
=\begin{pmatrix}\displaystyle\teta-\teta_c\\0\end{pmatrix}\quad\hbox{ a.e. in }Q_T,\\
\label{p4}
&\teta=\gamma(w),\quad\xib=(\xi_1,\xi_2)\in\partial\varphi(\betab)\quad\hbox{a.e.
in }Q_T,
\end{align}
and such that
\begin{align}
\label{p5}
&w(0)=w_0,\quad p(0)=p_0\quad\hbox{a.e. in }\Omega,\\
\label{p6}
&\betab(0)=(\beta_1(0),\beta_2(0))=\betab_0\quad\hbox{a.e. in }\Omega.
\end{align}

\beos\label{meanrho}
Notice that testing formally equation \eqref{p1} by $1$ one can immediately deduce
that the ``total mass'' of $p+\betab^\varrho$, where
$\betab^\varrho:=(\varrho_1\beta_1,\varrho_2\beta_2)$ is conserved  during the time interval $[0,T]$. Indeed, we have
\begin{equation}\label{meanbetaro}
\int_\Omega(p+\varrho_1\beta_1+\varrho_2\beta_2)(t)=\int_\Omega(p+\varrho_1\beta_1+\varrho_2\beta_2)(0)
\end{equation}
for all $t\in [0,T]$. The physical meaning of this relation is
that the mass of water either liquid or solid contained in the domain
$\Omega$ is constant, due to the watertight boundary condition we
have chosen (cf.~\eqref{bop} with $\alpha_p=0$). Finally, let us note that, by interpolation,
from \eqref{reg2} it follows that $w\in \CZH$ and so all the three boundary conditions in (\ref{p5}--\ref{p6})
hold in $H$.
\eddos

\beos\label{kappa}
Obviously, an admissible convex function satisfying assumption~\eqref{ipophi} is, e.g.,
$\varphi=I_{ K}$, being $K$ the convex set of $\RR^2$ defined in the
Introduction.
\eddos

\beos\label{tetab}
Let us underline here (cf.~also \cite[Remark~3.1]{bcf}) some regularities properties of
possible solutions to {\sc Problem (P)}. First of all, let us recall the
Gagliardo-Niremberg inequality in 3D, that is
\begin{equation}\label{gagli}
\Vert v\Vert_{L^p(\Omega)}\leq c\Vert v\Vert_V^\alpha\Vert v\Vert_H^{1-\alpha}\quad\forall v\in V,
\end{equation}
for $\frac{1}{p}=\frac{\alpha}{6}+\frac{1-\alpha}{2}$ and for some
positive constant $c$. Hence, if we are able to prove that
$\teta^{1/2}\in\LIH\cap\LDV$, then it follows that
\begin{equation}\no
\teta\in L^{5/3}(Q).
\end{equation}
\eddos

\beos\label{sigma}
Note that in our analysis we could also consider a more general form of our equation \eqref{p3};
we can  treat the case in which an antimonotone, but smooth,  contribution
is added to the monotone one given by $\partial\varphi$. Indeed, we can cover with our analysis
the case in which \eqref{p3} is substituted by the more general inclusion:
\begin{equation}
\displaystyle\dt{\betab}
+\nu\begin{pmatrix}{\cal B}\beta_1\\{\cal B}{\beta_2}\end{pmatrix}+\partial\varphi(\betab)+\sigma'(\betab)-
p\begin{pmatrix}\varrho_1\\\varrho_2\end{pmatrix}
\ni\begin{pmatrix}\displaystyle\teta-\teta_c\\0\end{pmatrix}\quad\hbox{ a.e. in }Q_T,
\end{equation}
where $\sigma\in C^{1,1}(\RR^2)$ is a, possibly antimonotone, contribution  coming from the
possibly non-convex part $\sigma(\betab)$ in the free energy functional \eqref{psi}.
The example we have in mind is the
one given by the ``classical'' double-well potential.
\eddos

The aim of the following Section~\ref{appro} is to prove one of our main results,
that is the following global existence theorem.

\bete\label{exiteo}
Let Hypothesis~\ref{hyp1} hold true and let $T$ be a positive final time. Then {\sc Problem (P)}
has at least a solution
on the whole time interval $[0,T]$.
\ente

\beos\label{unioss} The problem of finding also a uniqueness
result for {\sc Problem (P)} is still open and is strictly related
to the possibility of finding more regularity on the
$\teta$-component of the solution to {\sc Problem (P)}.  The reader can refer also to
\cite{bcf} and \cite{bf} for further comments on this topic in
case the voids are not admissible and only the two equations \eqref{p2} and \eqref{p3} (with $p=0$ and
$\varphi=I_C$, being $C=\{(\beta_1,\beta_2)\in [0,1]\,|\,\,\beta_1+\beta_2=1\}$)
are coupled.
\eddos

The proof of Theorem~\ref{exiteo} is based on the following scheme, which follows the idea of \cite{bcf}:
we first solve an approximating
problem in which $\gamma$ is  substituted  by a Lipschitz-continuous function,
then we make a priori estimates (independent of the approximation parameter), which will allow us to pass
to the limit by means of compactness and monotonicity arguments. In view of these considerations,
let us state here a preliminary result (whose proof will be given in Section~\ref{proofpreli}).

Consider the following assumptions on the data
\behy\label{hyp2} Let  (\ref{ipophi}--\ref{iponu}) and (\ref{datochi}) hold true and assume
moreover that:
\begin{align}\no
&\delta:=\gamma^{-1}\,:\,\RR\to \RR^2
\quad\hbox{is maximal monotone operator}\\
\label{gammalip}
&\hbox{with Lipschitz continuous inverse graph } \gamma,\\
\label{datowpiu}
&w_0\in H,\quad {\cal R}\in \LDVp,
\end{align}
where ${\cal R}$ is defined in (\ref{sorg3}).
\enhy

\bete\label{teolip}
Suppose now that Hypothesis~\ref{hyp2} holds true.
Then {\sc Problem (P)} has at least a solution $(w,\,\beta_1,\beta_2)$, $(\teta,\,\xi_1,\xi_2,\,p)$
satisfying the following regularity properties:
\begin{align}\label{regwpiu}
&w,\,\teta\in \CZH\cap\LDV,\\
\label{regbetapiu}
&\nu \beta_1,\,\nu\beta_2\in\LDW, \quad \xi_1,\,\xi_2\in\LDH.
\end{align}
Moreover, the components $w$, $\beta_1,\,\beta_2$, $p$, and $\teta$ of such a solution are  uniquely
determined.
\ente

\beos\label{posi} Observe that the main advantage of taking the
entropy balance equation \eqref{p2} instead of the internal energy balance
equation is that once one has solved the problem in some sense and
has found the temperature $\teta:=\gamma(w)$, it is automatically
positive because it stands in the image of the function $\gamma$
(cf.~\eqref{gammaexp}). Indeed in many cases it is difficult to
deduce this fact only from the internal energy balance equation
(cf., e.g., \cite{cspos} in order to see one example of these
difficulties). Let us note that within the small
 perturbations assumption the entropy balance and the classical heat equation are equivalent in mechanical terms
 (cf.~\cite{bcf, bf}).
\eddos

Moreover, in the same framework as above, let us recall the following inequality, holding in
$\Omega$ introduced above (cf.~\eqref{defiwo}):
\begin{align}
\label{poinc}
&\Vert v-v_\Omega\Vert_H\leq c\Vert\nabla v\Vert_H\quad \forall v\in V,
\end{align}
for some constants $c>0$, depending only on $\Omega$. Inequality \eqref{poinc}
 is one form of the standard Poincar\'e-Wirtinger inequality.
 Finally, we will make use of the following elementary inequality:
\begin{equation}\label{elem}
ab\leq\eta a^2+\frac{1}{4\eta}b^2,\quad\forall a,b\in\RR\quad\forall\eta>0.
\end{equation}

%%%%%%%%%%%%%%%%%%%%%%%%%%%%%%%%%%%%%%%%%%%%%%%%%%%%%%%%%%%%%%%%%%%%%%%%%%%%%%%%%%%
%%%%%%%%%%%%%%%%%%%%%%%%%%%%%%%%%% lipproof %%%%%%%%%%%%%%%%%%%%%%%%%%%%%%%%%%%%%%%

\section{Proof of Theorem~\ref{teolip}}
\label{proofpreli}

In this section we assume Hypothesis~\ref{hyp2} and we prove Theorem~\ref{teolip}.
Here we are going  to prove existence and uniqueness of solutions to Problem $(P)\ee$
in the time interval $[0,T]$, using a fixed point technique of contractive type.

In order to develop the proof in a rigorous way, we introduce
here the Lipschitz continuous Yosida-Moreau approximation
$\partial\varphi\ee=(\varphi\ee)'$
(cf.~\cite[Prop.~2.11, p.~39]{brezis}) of $\partial\varphi$, we first solve
the {\sc Problem (P)} with $\partial\varphi\ee$ instead of $\partial\varphi$ (call
it $(P)\ee$ in this section) and then we perform a
priori estimates independent of $\e$ and pass to the limit as
$\e\searrow 0$ recovering a solution to {\sc Problem (P)}. Moreover,
 we denote by $c$ the positive constants (maybe different from line to line), depending
on the data of the problem, but not on $\e$.
Then we are ready to introduce the variational
formulation of our approximating Problem $(P)\ee$ as follows.

\medskip
\noindent
{\sc Problem $(P)\ee$.} Find $(w\ee,\,\beta_1\ee,\beta_2\ee)$
and $(\teta\ee,\,\xi_1\ee,\xi_2\ee,\,p\ee)$ with the regularities
\begin{align}
\label{reg2e}
&w\ee,\,\teta\ee\in \HUVp\cap\LDV,\\
\label{reg3e}
&\betab\ee=(\beta_1\ee,\beta_2\ee)\in (\HUH)^2,
\quad \nu\beta_1\ee,\,\nu\beta_2\ee\in \LIV\cap L^2(0,T;W),\\
\label{reg4e}
&\xi_1\ee,\xi_2\ee\in L^{2}(Q_T),\quad p\ee\in\CZV\cap\LDW\cap\HUH,
\end{align}
satisfying
\begin{align}\label{p1e}
&\dt p\ee+\varrho_1\dt\beta_1\ee+\varrho_2\dt\beta_2\ee+{\cal B}p\ee=0\quad\hbox{ a.e. in }Q_T,\\
\label{p2e}
&\dt w\ee
+\dt\beta_1\ee+{\cal B}w\ee={\cal R}\quad\hbox{in }V',\hbox{ a.e. in }[0,T],\\
\label{p3e}
&\displaystyle\dt\betab\ee
+\nu\begin{pmatrix}{\cal B}\beta_1\ee\\{\cal B}{\beta_2\ee}\end{pmatrix}+\xib\ee
-p\ee\begin{pmatrix}\varrho_1\\\varrho_2\end{pmatrix}
=\begin{pmatrix}\displaystyle\teta\ee-\teta_c\\0\end{pmatrix}\qquad\hbox{ a.e. in }Q_T,\\
\label{p4e}
&\teta\ee=\gamma(w\ee),\quad\xib\ee=(\xi_1\ee,\xi_2\ee)=(\varphi\ee)'(\betab\ee)\quad\hbox{a.e. in }Q_T,
\end{align}
and such that
\begin{align}
\label{p5e}
&w\ee(0)=w_0,\quad p\ee(0)=p_0\quad\hbox{a.e. in }\Omega,\\
\label{p6e}
&\betab\ee(0)=(\beta_1\ee(0),\beta_2\ee(0))=\betab_0\quad\hbox{a.e. in }\Omega.
\end{align}

Now we start here proving existence of solutions to $(P)\ee$, by means
of a contraction argument.

\paragraph{Existence of solutions to $(P)\ee$, $\e\in (0,\bar\e]$.}
Let us take $\bar{t}\in [0,T]$
(we will choose it later) and call
$${\cal X}:=\{(r,s)\in (H^1(0,\bar{t};H))^2\,:\,(r,s)\in {\cal D}(\varphi)\}.$$
First fix $(\betapb,\betasb)\in {\cal X}$ in the equations \eqref{p1e} and \eqref{p2e}, then,
by well-known results (cf.~\cite[Lemma~6.3]{baiocchi} and also (\ref{sorg2})), we find also a unique
$w:={\cal T}_1(\betapb,\betasb)\in \HUH\cap\CZV\cap L^2(0,T;H^2(\Omega))$  solution of \eqref{p2e}
and $p:={\cal T}_2(\betapb,\betasb)\in \HUH\cap\CZV\cap L^2(0,T;W)$ solution to \eqref{p1e}.
Then, if we take this values of $w$  and $p$ in \eqref{p3e} in place of $w\ee$, $p\ee$,
again by standard results, we can find the solution $\betab=(\beta_1,\beta_2)$ to \eqref{p3e}.

In this way, we have defined an operator ${\cal T}:\,{\cal
X}\to{\cal X}$ such that $(\beta_1,\beta_2)=:{\cal
T}(\betapb,\betasb)$. What we have to do now is to prove that
${\cal T}$ is a contraction mapping on ${\cal X}$ for a
sufficiently small $\bar{t}\in [0,T]$. In order to prove that
${\cal T}$ is contractive, let us proceed by steps and forget of
the apices $\e$.
\smallskip

\noindent
{\sc First step.} Let $(\betapb^i,\betasb^i)\in{\cal X}$ ($i=1,2$),
$p^i={\cal T}_2(\betapb^i,\betasb^i)$,
$w^i={\cal T}_1(\betapb^i,\betasb^i)$,
and $(\beta_1^i,\beta_2^i)={\cal T}(\betapb^i,\betasb^i)$. Then,
writing two times \eqref{p2e} with $(\betapb^i,\betasb^i)$ ($i=1,2$), making the difference,
testing the resulting  equation with $(w^1-w^2)$, and integrating on $(0,t)$ with $t\in [0,\bar t]$,
we get the following inequality
\begin{align}\no
&\Vert w^1(t)-w^2(t)\Vert_H^2+\Vert w^1-w^2\Vert_{\LDtV}^2\\
\label{contr1}
&\leq C_1\Big(\Vert(\betapb^1)_t-(\betapb^2)_t\Vert_{\LDtH}^2+\Vert(\betasb^1)_t-(\betasb^2)_t\Vert_{\LDtH}^2\Big),
\end{align}
for some positive constant $C_1$ independent of $t$.
\smallskip

\noindent
{\sc Second step.} We write now  \eqref{p1e} with $(\betapb^i,\betasb^i)$,
make the difference, test the resulting  equation with $p^1-p^2$,
and integrate on $(0,t)$ with $t\in [0,\bar t]$, then we get the following inequality
\begin{align}\no
&\Vert (p^1-p^2)(t)\Vert_{H}^2+\Vert p^1-p^2\Vert_{L^2(0,t;V)}^2\\
\label{contr2}
&\leq C_2\Big(\Vert(\betapb^1)_t-(\betapb^2)_t\Vert_{\LDtH}^2+\Vert(\betasb^1)_t-(\betasb^2)_t\Vert_{\LDtH}^2\Big),
\end{align}
for some positive constant $C_2$ independent of $t$.

\smallskip
\noindent
{\sc Third step.} Let us take $\teta^i=\gamma(w^i)$ ($i=1,2$) and write equation \eqref{p3e} for $\teta^i$ and $u^i$,
make the difference between the two equations written for $i=1$ and $i=2$, and test the resulting vectorial equation
by the vector $\big((\beta_1^1)_t-(\beta_1^2)_t, (\beta_2^1)_t-(\beta_2^2)_t\big).$ Summing up the two lines and
integrating on $(0,t)$ with $t\in [0,\bar t]$, we have
\begin{align}\no
&\sum_{j=1}^2\Vert(\beta_j^1)_t-(\beta_j^2)_t\Vert_{\LDtH}^2+\nu\sum_{j=1}^2\Vert\nabla(\beta_j^1-\beta_j^2)(t)\Vert_H^2\\
\no
&\leq\int_{Q_T}\left|(p^1-p^2)\big(\varrho_1((\beta_1^1)_t-(\beta_1^2)_t)+\varrho_2((\beta_2^1)_t-(\beta_2^2)_t)\big)\right|\\
\no
&\quad+\frac14\sum_{j=1}^2\Vert(\beta_j^1)_t-(\beta_j^2)_t\Vert_{\LDtH}^2+2\Vert\teta^1-\teta^2\Vert_{\LDtH}^2\\
\no
&\leq \mezzo\sum_{j=1}^2\Vert(\beta_j^1)_t-(\beta_j^2)_t\Vert_{\LDtH}^2+2\Vert\teta^1-\teta^2\Vert_{\LDtH}^2
+2\Vert p^1-p^2\Vert_{\LDtH}^2.
\end{align}
Moreover, using the assumption \eqref{gammalip} on $\gamma$, we
get the following inequality
\begin{align}\no
&\mezzo\sum_{j=1}^2\Vert(\beta_j^1)_t-(\beta_j^2)_t\Vert_{\LDtH}^2
+\nu\sum_{j=1}^2\Vert\nabla(\beta_j^1-\beta_j^2)(t)\Vert_H^2\\
\label{contr3}
&\leq 2t\left(\Vert w^1- w^2\Vert_{\CZtH}^2+\Vert p^1- p^2\Vert_{\CZtH}^2\right).
\end{align}

\smallskip
\noindent
{\sc Fourth step.} Summing up the two inequalities \eqref{contr2} and \eqref{contr3},
and using \eqref{contr1}, we get
\begin{align}\no
&\sum_{j=1}^2\Vert(\beta_j^1)_t-(\beta_j^2)_t\Vert_{\LDtH}^2+\nu\sum_{j=1}^2\Vert\nabla(\beta_j^1
-\beta_j^2)(t)\Vert_H^2\\
\label{contr4}
&\leq C_3\,t\sum_{j=1}^2\Vert(\beta_j^1)_t-(\beta_j^2)_t\Vert_{\LDtH}^2,
\end{align}
for some positive constant $C_3$ independent of $t$. Hence, choosing $t$ sufficiently small (this is our ${\bar t}$),
we recover the contractive property of ${\cal T}$. Moreover,
applying the Banach fixed point theorem to ${\cal T}$, we get a
unique solution for the Problem $(P)\ee$ on the time interval
$[0,\bar t]$. Then, it is not difficult to check that
$$|{\cal T}(\betapb^1,\betasb^1)^k-{\cal T}(\betapb^2,\betasb^2)^k|_{(\HUH)^2}\leq\frac{(C_3T)^k}{k!}
|(\betapb^1,\betasb^1)-(\betapb^2,\betasb^2)|_{(\HUH)^2}$$
for any power ${\cal T}^k$ of ${\cal T}$. Hence, one can find $k$ such that $(C_3T)^k<k!$
and consequently the corresponding ${\cal T}^k$ is a contraction and admits a unique fixed point,
which is the searched solution to {\sc Problem $(P)\ee$.} on $[0,T]$.

\paragraph{Estimates and passage to the limit.}
We perform now uniform in $\e$ estimates on the solution to {\sc Problem $(P)\ee$}, which
allow us to pass to the limit as $\e\searrow 0$.
\medskip

\noindent
{\sc First estimate.}
Test \eqref{p1e} by $p$,
\eqref{p2e} by $\teta=\gamma(w)$,
and \eqref{p3e}  by $\dt\betab$. Then, sum up the result of the first two tests with the sum of the two
components of the vectorial equation found by the third test. This, after an integration in time on
$(0,t)$, with $t\in [0,T]$, thanks to a cancellation of two integrals, leads to
\begin{align}\no
&\frac12\Vert p(t)\Vert_{H}^2+\Vert\nabla p\Vert_{\LDtH}^2+\io\widehat\gamma(w(t))+
c_\gamma\itt\io|\nabla\teta|^2+\sum_{j=1}^2\Vert\dt\beta_j\Vert_{\LDtH}^2\\
\no
&+\frac{\nu}{2}\sum_{j=1}^2\Vert\nabla\beta_j(t)\Vert_{H}^2
+\io \varphi\ee(\betab(t))\leq \frac12\Vert p_0\|_H^2+\io\widehat\gamma(w_0)+\itt\duav{{\cal R},\teta}
\\
\label{sprimalip}
&+\frac{\nu}{2}\sum_{j=1}^2\Vert\nabla\beta_j(0)\Vert_H^2+ \io \varphi\ee(\betab_0)
\displaystyle+\itt \io\teta_c\dt\beta_1.
\end{align}
Now, using \eqref{datochi}, (\ref{gammalip}--\ref{datowpiu}), and Schwarz inequality with \eqref{elem},
  we get:
\begin{align}\no
&\Vert p(t)\Vert_{H}^2+\Vert\nabla p\Vert_{\LDtH}^2+\io\teta(t)+\itt\io|\nabla\teta(s)|^2\,ds
+\sum_{j=1}^2\Vert\dt\beta_j\Vert_{\LDtH}^2
\\
\no
&+\nu\sum_{j=1}^2\Vert\nabla\beta_j(t)\Vert_{H}^2
+\io \varphi\ee(\betab(t))\leq c+c_\eta\Vert R\Vert_{\LDtVp}^2+\eta\Vert\teta\ee\Vert_{\LDtV}^2,
\end{align}
holding true for all $\eta>0$ and for some positive $c_\eta$.
Finally, using \eqref{poinc} and \eqref{datowpiu}, and choosing $\eta$ sufficiently small, we get:
\begin{align}\no
&\Vert p(t)\Vert_{H}^2+\Vert\nabla p\Vert_{\LDtH}^2+
\Vert\teta\Vert_{\LDtV}^2+\sum_{j=1}^2\Vert\dt\beta_j\Vert_{\LDtH}^2\\
\label{s1lip}
&+\nu\sum_{j=1}^2\Vert\nabla\beta_j(t)\Vert_H^2
+\io \varphi\ee(\betab(t))\leq c.
\end{align}
Then, testing \eqref{p2e} with $w$ and using \eqref{s1lip}, it is
a standard matter to deduce the following bound
\begin{equation}\label{s2lip}
\Vert w\Vert_{\LIH\cap\LDV\cap\HUVp}\leq c,
\end{equation}
from which, using assumptions \eqref{gammalip} and \eqref{datochi}, we also deduce
\begin{equation}\label{s2lipbis}
\Vert \teta\Vert_{\LIH\cap\LDV\cap\HUVp}+\Vert p\Vert_{\HUH\cap\LIV\cap \LDW}\leq c.
\end{equation}

\noindent {\sc Second estimate.} Test now \eqref{p3e} by $\xib\ee$. Note that the term
$$
(\nabla\betab\ee(t),\nabla\xib\ee(t))=(\nabla\betab\ee(t),\nabla(\varphi\ee)'(\betab\ee(t)))
$$
is non-negative due to the monotonicity of $(\varphi\ee)'$.
Hence, due to estimate \eqref{s1lip} and \eqref{s2lipbis}, we obtain
\begin{equation}\label{s3lip}
\|\xi_i\ee\|_{\LDH}\leq c, \quad i=1,\,2,
\end{equation}
and also
\begin{equation}\label{s5lip}
\|\nu\beta_i\|_{\LDW}\leq c, \quad i=1,\,2,
\end{equation}
by comparison in \eqref{p3} and by standard regularity results
for elliptic equations.
\smallskip

\noindent
{\sc Passage to the limit.} As we have just mentioned, we want to conclude the proof
of Theorem~\ref{teolip} passing to the limit in Problem $(P)\ee$ as $\e\searrow 0$ using the previous uniform
(in $\e$) estimates on its solution and exploiting some compactness-monotonicity arguments.
Let us list before the weak or weak-star convergence coming directly from the previous estimates and well-known
weak-compactness results. Note that the following convergences hold only up to a subsequence of $\e$ tending to 0
(let us say $\e_k\searrow 0$). We denote it again with $\e$ only for simplicity of notation. From the estimates
(\ref{s1lip}--\ref{s2lipbis}), (\ref{s3lip}--\ref{s5lip}), we deduce that there exist
$(w,\bar\teta,p,\xi_j),$ $j=1,2$, such that
\begin{align}\label{c1lip}
w\ee\to w\quad&\hbox{weakly star in }H^1(0,T;V')\cap L^2(0,T; V),\\
\label{c2lip}
\teta\ee\to \bar\teta\quad&\hbox{weakly star in }H^1(0,T;V')\cap L^2(0,T; V),\\
\label{c3lip}
\beta_j\ee\to\beta_j\quad&\hbox{weakly in }H^1(0,T;H),\,(j=1,2),\\
\label{c4lip}
\nu\beta_j\ee\to\nu\beta_j\quad&\hbox{weakly star in }\LIV\cap L^{2}(0,T;W),\,(j=1,2),\\
\label{c5lip}
p\ee\to p\quad&\hbox{weakly star in }\HUH\cap \LIV\cap \LDW,\\
\label{c6lip}
\xi_j\ee\to \xi_j\quad&\hbox{weakly in }L^{2}(Q_T), \quad (j=1,2).
\end{align}
Moreover, employing the Aubin-Lions lemma (cf.~\cite[p.~58]{lions}) and
\cite[Cor.~8, p.89]{simon}, we also get
\begin{align}\label{c7lip}
w\ee\to w\quad&\hbox{strongly in }\LDH,\quad \hbox{and hence  a.e. in }Q_T,\\
\label{c7bis}
\teta\ee\to \bar\teta\quad&\hbox{strongly in }\LDH,\quad \hbox{and hence  a.e. in }Q_T,\\
\label{c8lip}
\beta_j\ee\to\beta_j\quad&\hbox{strongly in }\CZH\cap\LDV\quad\hbox{if }\nu>0\quad (j=1,2),\\
\no
p\ee\to p \quad &\hbox{strongly in }\CZH\cap \LDV\quad \hbox{and hence  a.e. in }Q_T.
\end{align}
Note that convergences (\ref{c2lip}), \eqref{c7lip}, and relation \eqref{pe4} imply immediately
(cf.~\cite[Lemme~1.3, p.~12]{lions})
the convergence
\begin{equation}\no
\teta\ee\to\bar\teta=\gamma(w)\quad\hbox{a.e. in }Q_T.
\end{equation}
Finally, it remains to prove the identification of the maximal monotone graph $\partial\varphi$, i.e.
\begin{equation}\no
\xib\in\partial\varphi(\betab)\quad\hbox{a.e. in }Q_T,
\end{equation}
with $\xib=(\xi_1,\xi_2)$ and $\xi_j$ $(j=1,2)$ that are the weak limits defined in \eqref{c6lip}.
In order to do that we
should verify that
\begin{equation}\label{limsuplip}
\limsup_{\e\searrow 0}\int_0^T(\xib\ee,\betab\ee)\leq \int_0^T(\xib,\betab).
\end{equation}
The proof here is splitted into two parts.

\paragraph{Case $\nu>0$.} From \eqref{c8lip} and \eqref{c6lip} we immediately deduce that \eqref{limsup}
is verified.

\paragraph{Case $\nu=0$.}  In this
case we can prove that $\betab\ee$ is a Cauchy sequence in $\CZH$.
Indeed one can take the differences of equations \eqref{p3}
written for two different indices $\e$ and $\e'$ and test it by
the difference vector $(\betab\ee-\betab^{\e'})$. Take the
differences of the two equations \eqref{p1} integrated in time and test it by $p^{\e}-p^{\e'}$.
Summing up the two resulting equations and integrating over
$(0,t)$,  $t\in (0,T)$, we get:
\begin{align}\no
&\|p\ee-p^{\e'}\|_{\LDtH}^2+\mezzo\|(\betab\ee-\betab^{\e'})(t)\|_{H}^2
\leq\int_{Q_t}(\teta\ee-\teta^{\e'})(\beta_1\ee-\beta_1^{\e'})\\
\no
&\leq \itt\|(\teta\ee-\teta^{\e'})(s)\|_{H}
\|(\beta_1\ee-\beta_1^{\e'})(s)\|_{H}\,ds.
\end{align}
Applying now the Gronwall lemma \cite[Lemme~A.3]{brezis}, we get
\begin{align}\no
&\|(\betab\ee-\betab^{\e'})(t)\|_{H}^2\leq c\|\teta\ee-\teta^{\e'}\|_{L^1(0,T;H)}\to 0\quad\hbox{as }
\e,\,\e'\searrow 0.
\end{align}
Hence, we have that $\betab\ee\to \betab$ strongly in $\CZH$, which is sufficient in order to
prove that \eqref{limsup} is satisfied.

Hence, in both cases, all these convergences with the
identifications made above make us able to pass to the limit (as
$\e\searrow 0$ or at least for a subsequence of it) in Problem
$(P)\ee$ finding a solution to {\sc Problem (P)}. Note that once we prove uniqueness,
the convergences above will turn out to hold for all the sequence $\e\searrow 0$ and not only up to a
subsequence.

\paragraph{Uniqueness.} Assume that $(w^1,\,\betab^1), (\teta^1,\,\xib^1,\, p^1)$ and
$(w^2,\,\betab^2), (\teta^2,\,\xib^2,\, p^2)$  are two solutions of {\sc Problem (P)} corresponding to
the same data, and let us call by $\overline{w}:=w^1-w^2,\,\overline p:=p^1-p^2,\,\overline{\beta}_j:=\beta_j^1-\beta_j^2,\, (j=1,2)$.
Then make the differences from \eqref{p1} written down for the first solution and
\eqref{p1} integrated in time  and written down for the second solution, then test the resulting equation by
$\overline p$, test the difference between the equations \eqref{p3} written in terms of $\betab^1$ and $\betab^2$ by
$(\overline\beta_1,\overline\beta_2)$. Finally,  integrate in time \eqref{p2}
(let us call it $1*\eqref{p2}$) and make the differences from $1*\eqref{p2}$ written down for the first solution and
$1*\eqref{p2}$ written down for the second solution, then test the resulting equation by
$\overline w$. Then sum up the three resulting equations and integrate it over $(0,t)$, use
the maximal monotonicity of $\partial\fhi$, getting
\begin{align}\no
&\|\overline p\|_{\LDtH}^2+\left\|\itt \nabla\overline{p}\right\|_{H}^2+\sum_{j=1}^2\|\overline{\beta}_j(t)\|_{H}^2
+\nu\sum_{j=1}^2\|\nabla\overline{\beta}_j\|_{\LDtH}^2\\
\no
&+\|\overline{w}\|_{\LDtH}^2
+\left\|\itt \nabla\overline{w}\right\|_{H}^2\leq c\itt(\gamma(w^1)-\gamma(w^2)-\overline{w},\overline{\beta}_1).
\end{align}
Then, getting advantage
of the Lipschitz continuity of $\gamma $ (cf.~assumption~\eqref{gammalip}) and using \eqref{elem},
we get
\begin{align}\no
&\|\overline p\|_{\LDtH}^2+\left\|\itt \nabla\overline{p}\right\|_{H}^2+\sum_{j=1}^2\|\overline{\beta}_j(t)\|_{H}^2
+\nu\sum_{j=1}^2\|\nabla\overline{\beta}_j\|_{\LDtH}^2\\
\no
&+\|\overline{w}\|_{\LDtH}^2
+\left\|\itt \nabla\overline{w}\right\|_{H}^2\leq \eta \|\overline{w}\|_{\LDtH}^2
+c_\eta\itt\|\overline{\beta}_1\|_H^2,
\end{align}
holding true for all positive constants $\eta$ and for some $c_\eta>0$. Finally, choosing
$\eta$ sufficiently small and applying a standard version of
the Gronwall lemma (cf.~\cite [Lemme~A.2]{brezis}), we get the desired uniqueness of
the $\betab$ and $w$-components of the solution to {\sc Problem (P)} and
this concludes the proof of Theorem~\ref{teolip}.\fin

%%%%%%%%%%%%%%%%%%%%%%%%%%%%%%%%%%%%%%%%%%%%%%%%%%%%%%%%%%%%%%%%%%%
%%%%%%%%%%%%%%%%%%%%%%%%%%%% exiproof  %%%%%%%%%%%%%%%%%%%%%%%%%%%%%%%%

\section{Proof of Theorem~\ref{exiteo}}
\label{exiproof}

The following section is devoted to the proof of
Theorem~\ref{exiteo}. It will be done by steps. First we
approximate our {\sc Problem (P)} by a more regular Problem $(P)\ss$,
then (fixed $\s>0$) we find well-posedness for the approximating
problem using a fixed-point theorem and then we perform some
a-priori estimates (independent of $\s$) on its solution. Finally
we find existence of solutions of {\sc Problem (P)} by passing to the
limit  in Problem $(P)\ss$ as $\s\searrow 0$.

In order to perform the following estimates rigorously, we use here
the Lipschitz continuous Yosida-Moreau approximation
$\partial\varphi\ss=(\varphi\ss)'$
(cf.~\cite[Prop.~2.11, p.~39]{brezis}) of $\partial\varphi$.

\subsection{The approximating problem}\label{appro}

Let us introduce here the approximating Problem $(P)\ss$ of our {\sc Problem (P)} (cf.~equations
(\ref{p1}--\ref{p6})).
First of all,
we take a small positive parameter $\s$ and we call
$\gamma\ss$ the following Lipschitz continuous approximation of the function
$\gamma(w)=\exp(w)$ in \eqref{gammaexp}, i.e.~the function
\begin{equation}\label{gamma}
\gamma\ss(w):=\begin{cases}
                    \exp w&\hbox{if }r\leq 1/\s\\
                    (w-1/\s)\exp(1/\s)+\exp(1/\s)&\hbox{if }r\geq 1/\s.
                 \end{cases}
\end{equation}
Moreover let $\delta\ss$ be the inverse function of $\gamma\ss$, i.e.
\begin{equation}\label{delta}
\delta\ss\big(\gamma\ss(r)\big)=r\quad\forall r\in\RR
\end{equation}
and let $\widehat\gamma\ss$ be a primitive of the function $\gamma\ss$, i.e.
\begin{equation}\label{gammaprim}
\widehat\gamma\ss(r)=1+\int_0^r\gamma\ss(s)\,ds\quad\forall r\in \RR.
\end{equation}
Then the following properties of $\gamma\ss$ hold true (cf.~also \cite[Lemma~5.1]{bcf}).
\bele\label{progamma}
There holds
$$\widehat\gamma\ss(r)\geq\gamma\ss(r),\quad r(\delta\ss)'(r)\geq 1\quad\forall r\in\RR.$$
\enle
Now, we can approximate {\sc Problem (P)} as follows.
\smallskip

\noindent
{\sc Problem $(P)\ss$.} Find $(w\ss,\,\beta_1\ss,\beta_2\ss)$ and $(\teta\ss,\,\xib\ss,\,p\ss)$ with the following
regularity properties:
\begin{align}
\label{rege2}
&w\ss,\,\teta\ss\in\HUVp\cap\LDV,\\
\label{rege3}
&\beta_1\ss,\,\,\beta_2\ss\in \HUH,\quad\nu\beta_1\ss,\,\nu\beta_2\ss\in\LIV\cap\LDW, \\
\label{rege4}
&\xib\ss\in \LDH,\quad p\ss\in \HUH\cap\LIV\cap\LDW,
\end{align}
satisfying
\begin{align}\label{pe1}
&\dt p\ss+\varrho_1\dt\beta_1\ss+\varrho_2\dt\beta_2\ss+{\cal B}p\ss=0\quad\hbox{ a.e. in }Q_T,\\
\label{pe2} &\dt w\ss
+\dt\beta_1\ss+{\cal B} w\ss={\cal R}\quad\hbox{in }V'\hbox{ and a.e. in }[0,T],\\
\label{pe3}
&\dt{\betab\ss}+\nu
\begin{pmatrix}
{\cal B}{\beta_1}\ss\\{\cal B}{\beta_2}\ss\end{pmatrix}+\xib\ss
-p\ss\begin{pmatrix}\varrho_1\\\varrho_2
\end{pmatrix}
=
\begin{pmatrix}
\teta\ss-\teta_c\\0
\end{pmatrix}
\quad\hbox{a.e. in }Q_T,\\
\label{pe4}
&\teta\ss=\gamma\ss(w\ss),\quad\xib\ss\in \partial\varphi\ss(\betab\ss),\quad\hbox{a.e. in }Q_T,
\end{align}
and such that
\begin{align}\label{pe6}
&w\ss(0)=w_0,\quad p\ss(0)=p_0\quad\hbox{a.e. in }\Omega,\\
\label{pe7}
&\betab\ss(0)=(\beta_1\ss(0),\beta_2\ss(0))=\betab_0\quad\hbox{a.e. in }\Omega.
\end{align}

Existence and uniqueness of solutions to {\sc Problem $(P)\ss$} directly follow from Theorem~\ref{teolip}.
We proceed now performing a-priori estimates on this solution uniform in $\s$ in order to pass to
the limit in  {\sc Problem $(P)\ss$}  as $\s\searrow 0$, recovering a solution to
{\sc Problem $(P)$}.

\subsection{A priori estimates}\label{esti}

In this subsection we perform a-priori estimates on Problem $(P)\ss$ uniformly in $\s$, which will lead us
pass to the limit in Problem $(P)\ss$ as $\s\searrow 0$ and recover a solution of {\sc Problem (P)}. Hence,
let us denote by $c$ all the positive constants (which may also differ from line to line) independent
of $\s$ and depending on the data of the problem.
\noindent
Now, let us come to the (uniform in $\s$) estimates on the solution to Problem $(P)\ss$.

\smallskip

\noindent
{\sc First estimate.} Test \eqref{pe1} by $p\ss$, \eqref{pe2} by $\teta\ss=\gamma\ss(w\ss)$,
and \eqref{pe3}  by $\dt\betab\ss$. Then, sum up the result of the first two tests with the sum of the two
components of the vectorial equation found by the third test. This, after an integration in time on
$(0,t)$, with $t\in [0,T]$, thanks to a cancellation of two integrals, leads to
\begin{align}\no
&\frac12\Vert p\ss(t)\Vert_{H}^2+\Vert\nabla p\ss\Vert_{\LDtH}^2+\io\widehat\gamma\ss(w\ss(t))+
\itt\io\nabla\delta\ss(\teta\ss)\nabla\teta\ss\\
\no
&+\sum_{j=1}^2\Vert\dt\beta_j\ss\Vert_{\LDtH}^2+\frac{\nu}{2}\sum_{j=1}^2\Vert\nabla\beta_j\ss(t)\Vert_{H}^2
+\io \varphi\ss(\betab\ss(t))\leq \frac12\|p_0\|_H^2+\io\widehat\gamma\ss(w_0)
\\
\label{sprima}
&+\itt\duav{{\cal R},\teta\ss}+\frac{\nu}{2}\sum_{j=1}^2\Vert\nabla\beta_j\ss(0)\Vert_H^2+ \io \varphi\ss(\betab_0)
\displaystyle +\itt\io \teta_c\dt\beta_1^\sigma.
\end{align}
Now, following the line of \cite[(5.5)--(5.7), p.~1583]{bcf}, we can deal with the source term ${\cal R}$ recalling
\eqref{sorg3} and using a well-known compactness inequality (cf.~\cite[Theorem~16.4]{lima}) in this way
\begin{align}\label{erre}
\itt\io R\teta\ss&\leq\itt\Vert R(s)\Vert_{L^{\infty}(\Omega)}\Vert\teta\ss(s)\Vert_{L^1(\Omega)}\,ds,\\
\no
\itt\int_{\partial\Omega} \Pi\teta\ss&\leq
c\Vert \Pi\Vert_{L^{\infty}(\Sigma)}\Vert (\teta\ss)^{1/2}\Vert_{L^2(\Sigma)}\\
\label{acca}
&\leq \zeta\Vert\nabla(\teta\ss)^{1/2}\Vert_{\LDtH}^2+C_\zeta\Vert(\teta\ss)^{1/2}\Vert_{\LDtH}^2
\end{align}
for all $\zeta>0$ and for some positive $c_\zeta$ depending also on $\Vert \Pi\Vert_{L^\infty(\Omega)}$
(cf.~\eqref{sorg2}) and on $\Omega$.
Now, collecting estimates (\ref{erre}--\ref{acca}), with $\zeta=1$, using Lemma~\ref{progamma}
with assumptions  (\ref{datoteta}--\ref{sorg2}), we get the inequality:
\begin{align}\no
&\frac12\Vert p\ss(t)\Vert_{H}^2+\Vert\nabla p\ss\Vert_{\LDtH}^2+2\io\teta\ss(t)
+\itt\io\frac{|\nabla\teta\ss(s)|^2}{\teta\ss(s)}\,ds
\\
\no
&+\frac{1}{2}\sum_{j=1}^2\Vert\dt\beta_j\ss\Vert_{\LDtH}^2+\frac{\nu}{2}\sum_{j=1}^2\Vert\nabla\beta_j\ss(t)\Vert_{H}^2
+\io \varphi\ss(\betab\ss(t))\\
\no
&\leq c+\itt\Vert R(s)\Vert_{L^{\infty}(\Omega)}\Vert\teta\ss(s)\Vert_{L^1(\Omega)}\,ds
+c\Vert(\teta\ss)^{1/2}\Vert_{\LDtH}^2
.
\end{align}
Finally, a standard version of Gronwall lemma gives:
\begin{align}\no
&\Vert p\ss\Vert_{\LItH}^2+\Vert\nabla p\ss\Vert_{\LDtH}^2+\Vert(\teta\ss)^{1/2}(t)\Vert_H^2+
\itt\Vert\nabla(\teta\ss)^{1/2}\Vert_{H}^2+\sum_{j=1}^2\Vert\dt\beta_j\ss\Vert_{\LDtH}^2\\
\label{s1}
&+\sum_{j=1}^2\Vert\nabla\beta_j\ss(t)\Vert_H^2
+\io \varphi\ss(\betab\ss(t))\leq c.
\end{align}

\noindent
{\sc Second estimate.} Testing \eqref{pe2} with $w\ss$ and using \eqref{s1}, it is
a standard matter to deduce the following bound
\begin{equation}\label{s2}
\Vert w\ss\Vert_{\LIH\cap\LDV}\leq c,
\end{equation}
whence \eqref{pe2} entails
\begin{equation}\label{s3}
\| w\ss\|_{\HUVp}\leq c,
\end{equation}
while, recalling  \eqref{gagli}, from \eqref{s1} we also get
\begin{equation}\label{s4}
\|\teta\ss\|_{L^{5/3}(Q_T)}\leq c.
\end{equation}
Moreover, testing equation \eqref{p3e} by $|\xib\ss|^{-1/3}\xib\ss$, we get
\begin{equation}\label{s6}
\Vert\xi_i\ss\Vert_{L^{5/3}(Q_T)}\leq c\quad (i=1,2),
\end{equation}
and so, by comparison in \eqref{pe3}, we get
\begin{equation}\label{s7}
\Vert\nu\beta_i\ss\Vert_{L^{5/3}(0,T;W^{2,5/3}(\Omega))}\leq c\quad (i=1,2).
\end{equation}
Now it remains only to pass to the limit in (\ref{pe1}--\ref{pe7}) as $\sigma\searrow 0$.
This will be the aim of the next subsection.

\subsection{Passage to the limit}\label{conclu}

As we have just mentioned, we want to conclude the proof of Theorem~\ref{exiteo} passing to the limit in the
well-posed (cf.~Subsection~\ref{appro}) Problem $(P)\ss$ as $\s\searrow 0$ using the previous uniform
(in $\s$) estimates on its solution (cf.~Subsection~\ref{esti}) and exploiting some compactness-monotonicity argument.
Let us list before the weak or weak-star convergence coming directly from the previous estimates and well-known
weak-compactness results. Note that the following convergences hold only up to a subsequence of $\s\searrow 0$
(let us say $\s_k\searrow 0$). We denote it again with $\s$ only for simplicity of notation. From the estimates
(\ref{s1}--\ref{s7}), we deduce that
\begin{align}\label{c1}
w\ss\to w\quad&\hbox{weakly star in }H^1(0,T;V')\cap L^2(0,T; V),\\
\label{c2}
\teta\ss\to\teta\quad&\hbox{weakly in }L^{5/3}(Q_T),\\
\label{c3}
\beta_j\ss\to\beta_j\quad&\hbox{weakly in }H^1(0,T;H),\,(j=1,2),\\
\label{c4}
\nu\beta_j\ss\to\nu\beta_j\quad&\hbox{weakly star in }\LIV\cap L^{5/3}(0,T;W^{2,5/3}(\Omega)),\,(j=1,2),\\
\label{c5}
p\ss\to p\quad&\hbox{weakly star in }\HUH\cap\LIV\cap\LDW,\\
\label{c6}
\xi_j\ss\to \xi_j\quad&\hbox{weakly in }L^{5/3}(Q_T), \quad (j=1,2).
\end{align}
Moreover, employing the Aubin-Lions lemma (cf.~\cite[p.~58]{lions}) and
\cite[Cor.~8, p.89]{simon}, we also get:
\begin{align}\label{c7}
w\ss\to w\quad&\hbox{strongly in }\LDH,\quad \hbox{and hence  a.e. in }Q_T,\\
\label{c8}
\beta_j\ss\to\beta_j\quad&\hbox{strongly in }\CZH\cap\LDV\quad\hbox{if }\nu>0\quad (j=1,2),\\
\label{c8bis}
p\ss\to p \quad&\hbox{strongly in }\CZH\cap\LDV\quad\hbox{and hence  a.e. in }Q_T.
\end{align}
Note that (\ref{c2}), \eqref{c7}, and relation \eqref{pe4} imply immediately (cf.~\cite[Lemme~1.3, p.~12]{lions})
the convergence
\begin{equation}\no
\teta\ss\to\teta=\gamma(w)\quad\hbox{a.e. in }Q_T,
\end{equation}
and by Egorov theorem, we also deduce
\begin{equation}\label{c9}
\teta\ss\to\teta\quad\hbox{strongly in }L^{q}(Q_T)\quad\forall q\in [1,5/3).
\end{equation}
Finally, it remains to prove the identification of the maximal monotone graph $\partial\varphi$, i.e.
\begin{equation}\no
\xib\in\partial\varphi(\betab)\quad\hbox{a.e. in }Q_T,
\end{equation}
with $\xib=(\xi_1,\xi_2)$ and $\xi_j$ $(j=1,2)$ that are the weak limits defined in \eqref{c6}.
In order to do that we
should verify that
\begin{equation}\label{limsup}
\limsup_{\s\searrow 0}\int_0^T(\xib\ss,\betab\ss)\leq \int_0^T(\xib,\betab).
\end{equation}
The proof here is splited into two parts.

\paragraph{Case $\nu>0$.} In this case from (\ref{c3}--\ref{c4}), using again
\cite[Cor.~8, p.~89]{simon}, we also deduce
\begin{equation}\label{c10}
\beta_j\ss\to\beta_j\quad\hbox{strongly in }L^r(0,T;L^q(\Omega)),\quad\forall r\in [1,+\infty),
\, q\in [1,6), \,(j=1,2).
\end{equation}
>From this convergence and from \eqref{c6} we immediately deduce that \eqref{limsup}
is verified.

\paragraph{Case $\nu=0$.} In case $\nu=0$, due to assumption \eqref{iponu},
we also have that $\betab$ is bounded in $L^\infty(Q_T)$. In this
case we can prove that $\betab\ss$ is a Cauchy sequence in $\CZH$.
Indeed one can take the differences of equations \eqref{pe3}
written for two different indices $\s$ and $\s'$ and test it by
the difference vector $\betab\ss-\betab^{\s'}$. Take the
differences of the two equations \eqref{pe1} integrated in time and test it by $p\ss-p^{\s'}$.
Summing up the two resulting equations and integrating over
$(0,t)$,  $t\in (0,T)$, we get (as $\s,\,\s'\searrow 0$):
\begin{align}\no
\|p\ss-p^{\s'}\|_{\LDtH}^2+\mezzo\|(\betab\ss-\betab^{\s'})(t)\|_{H}^2
&\leq\int_{Q_t}(\teta\ss-\teta^{\s'})(\beta_1\ss-\beta_1^{\s'})\\
\no
&\leq 2\|\teta\ss-\teta^{\s'}\|_{L^1(Q_t)}\to 0.
\end{align}
Hence, we have that $\betab\ss\to \betab$ strongly in $\CZH$ and weakly star in $L^\infty(Q_T)$ and so
also strongly in $L^r(Q_T)$ for all $r\in [1,+\infty)$, which is sufficient in order to
prove that \eqref{limsup} is satisfied.

Hence, in both cases, all these convergences with the
identifications made above make us able to pass to the limit (as
$\s\searrow 0$ or at least for a subsequence of it) in Problem
$(P)\ss$ finding a solution to {\sc Problem (P)} and concluding in
this way the proof of Theorem~\ref{exiteo}. Note that the
convergences hold only for proper subsequences $\s_k$ of $\s$
tending to 0 because of lack of uniqueness of solutions (cf.~also
Remark~\ref{unioss}). This concludes to proof of
Theorem~\ref{exiteo}.\fin

%%%%%%%%%%%%%%%%%%%%%%%%%%%%%%%%%%%%%%%%%%%%%%%%%%%%%%%%%%%%%%%%%%%%%%%%%%%%%%%%%%%
%%%%%%%%%%%%%%%%%%%%%%%%%% biblio %%%%%%%%%%%%%%%%%%%%%%%%%%%%%%%%%%%%%%%%%%%%%%%

%%%%%%%%%%%%%%%%%%%%%%%%%%%%%%%%%%%%%%%%%%%%%%%%%%%%%%%%%%%%%%%%%%%
%%%%%%%%%%%%%%%%%%%%%%%%%%%%%%%%%%%%%%%%%%%%%%%%%%%%%%%%%%%%%%%%%%%

\end{document}